\documentclass[12pt]{amsart}
\title[soluciones discretas para sistemas matriciales en derivadas
parciales...]{soluciones discretas para sistemas matriciales en
derivadas parciales hiperb\'{o}licos y singulares}
\author{Manuel J. Salazar, Edison E. Villa}
\address{Universidad de Antioquia, Facultad de Ciencias Exactas y Naturales.
Departamento de Matematicas, A.A. 1226 Medellin, Colombia.}
\email{mjsalazar@matematicas.udea.edu.co,
eevc03@matematicas.udea.edu.co}
\date{\today }
\subjclass[2000]{35A09, 35L05}
\usepackage{amsfonts}
\usepackage{amsmath}
\usepackage{amssymb}
\usepackage[english,spanish]{babel}
\usepackage{graphicx,fourier}
\usepackage{pstricks,pst-node,pst-tree}

\newtheorem{theorem}{Teorema}

\newtheorem{definition}[theorem]{Definición}

\newtheorem{proposition}[theorem]{Proposición}

\begin{document}\maketitle

\selectlanguage{english}
\begin{abstract}
In this paper we study the construction of a discrete solution for a
hyperbolic system of partial differentials of the strongly coupled
type. In its construction, the discrete separation of matricial
variable method was followed. Two separate equations in differences
were obtained: a singular matricial and the other one a Sturm
Liouville vectorial problem, which by the superposition principle
yield a stable discrete solution.
\newline
\newline
Key Words: Hyperbolic system, boundary and Sturm Liouville problem,
mixed problem, stable solution.

\end{abstract}

\thispagestyle{empty} \selectlanguage{english}

\selectlanguage{spanish}

\begin{abstract}
En este trabajo se construye una soluci\'{o}n discreta a un sistema
en derivadas parciales hiperb\'{o}lico de tipo singular fuertemente
acoplado. En su construcci\'{o}n se sigui\'{o} el m\'{e}todo de
separaci\'{o}n de variables matricial discreto, obteni\'{e}ndose dos
ecuaciones en diferencias por separado: una matricial singular y la
otra un problema de Sturm Liouville vectorial, las cuales mediante
el principio de superposici\'{o}n, producen una soluci\'{o}n
discreta estable del sistema.
\newline
\newline
Palabras claves: Sistema hiperb\'{o}lico, contorno discretizado,
problema mixto y de Sturm Liouville, soluci\'{o}n estable.
\end{abstract}

\section{INTRODUCCION}

Los sistemas hiperb\'{o}licos en derivadas parciales con coeficientes
matriciales aparecen a menudo en la modelaci\'{o}n de situaciones tales como
el estudio de los de procesos de calefacci\'{o}n en microondas \cite{s2, S1}%
, \'{o}ptica \cite{S3}, cardiolog\'{\i}a \cite{S7a}, fen\'{o}menos s\'{\i}%
smicos en medios el\'{a}sticos \cite{S11}, entre otros. En este art\'{\i}%
culo consideraremos sistemas hiperb\'{o}licos en derivadas parciales con
coeficientes matriciales constantes, definidos en el rectangulo $[0,1]\times
\lbrack 0,T]$, los cuales caracterizan la ecuaci\'{o}n de onda homog\'{e}nea:%
\begin{equation}
Eu_{tt}(x,t)-Au_{xx}(x,t)=0,\ x\in \lbrack 0,1],\text{ }0<t<T,  \label{a}
\end{equation}%
\begin{equation}
\begin{array}{c}
A_{1}u(0,t)+A_{2}u_{x}(0,t)=0,\ 0<t<T, \\
B_{1}u(1,t)+B_{2}u_{x}(1,t)=0<t<T,%
\end{array}
\label{b}
\end{equation}%
\begin{equation}
\begin{array}{c}
u(x,0)=f(x),\ x\in \lbrack 0,1],\text{\ } \\
u_{t}(x,0)=g(x),\text{ }x\in \lbrack 0,1],%
\end{array}%
\text{\ }  \label{c}
\end{equation}%
donde $u=(u_{1},...,u_{m})^{T},$ $u_{tt}=(u_{1tt},...,u_{mtt})^{T},$ $%
u_{xx}=(u_{1xx},...,u_{mxx})^{T};$ $%
f=(f_{1},...,f_{m})^{T},g=(g_{1},...,g_{m})^{T}\in \mathbb{C}^{m},$ $A,$ $%
A_{1}$, $A_{2},$ $B_{1},$ $B_{2}$, $E$ $\in \mathbb{C}$ $^{m\times m},$ en
su estudio $\left( \text{\ref{b}}\right) $ es llamada condici\'{o}n de
contorno fuertemente acoplada y $\left( \text{\ref{c}}\right) $ el problema
mixto asociado a $\left( \text{\ref{a}}\right) -\left( \text{\ref{b}}\right)
$. El sistema $\left( \text{\ref{a}}\right) -\left( \text{\ref{c}}\right) $
es tambi\'{e}n llamado singular si $E$ es singular de lo contraio es
conciderado no singular. Una gama de casos no singulares de $\left( \text{%
\ref{a}}\right) -\left( \text{\ref{c}}\right) $ han sido ampliamente
tratados (v\'{e}ase \cite{P34}, \cite{camacho}, \cite{25a}, \cite{P34}, \cite%
{S10}). Si bien la literatura sobre su caso singular es muy poca \cite{PON}
estudia una mejora a \cite{P37a} la cual soluciona un caso particular de $%
\left( \text{\ref{a}}\right) -\left( \text{\ref{c}}\right) $. Este art\'{\i}%
culo se propone formalizar y generalizar los resultados de \cite{PON}
analizando la flexibilidad de las hip\'{o}tesis que se asumieron para la
obtenci\'{o}n de una soluci\'{o}n al sistema de ecuaciones mencionado. El
orden a desarrollar en el presente art\'{\i}culo es el siguiente:

-Secci\'{o}n 2: Se presentan definiciones y teoremas generales de las teor%
\'{\i}a de sistemas matriciales singulares y de las ecuaciones en
diferencias finitas, en particular demostraremos $3$ proposiciones
indispensables para la soluci\'{o}n de la versi\'{o}n discreta de $\left(
\text{\ref{a}}\right) -\left( \text{\ref{b}}\right) .$

-Secci\'{o}n 3: Luego de discretizar el sistema$\left( \text{\ref{a}}\right)
-\left( \text{\ref{c}}\right) $, aplicando el m\'{e}todo de separaci\'{o}n
de variables matricial discreto, se construye una soluci\'{o}n a la ecuaci%
\'{o}n de onda considerando su condici\'{o}n de contorno y se enuncia el
teorema que resume el resultado obtenido.

-Secci\'{o}n 4: Se estudia el problema mixto y la estabilidad de la soluci%
\'{o}n general. Con base en el teorema anterior enunciamos el segundo que
garantiza la existencia de una soluci\'{o}n discreta y estable de $\left(
\text{\ref{a}}\right) -\left( \text{\ref{c}}\right) $.

-Finalmente complementaremos lo obtenido con un ejemplo ilustrativo y
enunciaremos algunas observaciones y conclusi\'{o}n finales.

\section{PRELIMINARES}

En esta secci\'{o}n damos algunos resultados y notaciones que ser\'{a}n de
utilidad para las siguientes

\begin{definition}
\label{Def1}Sea A$\in
\mathbb{C}
^{m\times m},$ una matriz, la cual denotamos por $A^{G},$ es llamada inversa
generalizada de A si verifica la condici\'{o}n%
\begin{equation*}
AA^{G}A=A
\end{equation*}
\end{definition}

\begin{theorem}
\label{R1.2} $\left( Teorema\text{ }de\text{ }Mitra\right) $ Sea $A\in
\mathbb{C}
^{m\times m}$ y $A^{G}$ una inversa generalizada para A. Sea $b\in
\mathbb{C}
^{m}$, entonces%
\begin{equation*}
Ax=b\text{ tiene soluci\'{o}n}\Leftrightarrow AA^{G}b=b
\end{equation*}%
y cualquiera de sus soluciones tiene la forma
\begin{equation}
x=A^{G}b+\left( I-A^{G}A\right) z,\ \ \ z\in
\mathbb{C}
^{m}.  \label{ss7}
\end{equation}%
Demostraci\'{o}n: V\'{e}ase \cite{P340a}.
\end{theorem}

\begin{definition}
Sea W un subespacio de $%
\mathbb{C}
^{m}$ y A una matriz en $%
\mathbb{C}
^{m\times m},$ se dice que W es invariante por la matriz A, si verifica que $%
AW\subset W$ $\left( i.e\text{ }\forall x\in W\Rightarrow Ax\in W\right) .$
\end{definition}

\begin{theorem}
\label{T4} Sean $A,$ $B,$ $B^{G}\in
\mathbb{C}
^{m\times m},$ entonces $BA\left( I-B^{G}B\right) =0$ sii $Ker\left(
B\right) $ es un subespacio invariante de $A.$\newline
Demostraci\'{o}n: V\'{e}ase \cite{P340a}.
\end{theorem}

\begin{definition}
(La Inversa Drazin) Sea $A\in
\mathbb{C}
^{m\times m}$ tal que $ind(A)=k$, la inversa Drazin de A ($A^{D}$) es la
\'{u}nica matriz que cumple\newline
1) $A^{D}AA^{D}=A^{D}$\newline
2) $A^{D}A=AA^{D}$\newline
3) $A^{k+1}A^{D}=A^{k}$
\end{definition}

\begin{theorem}
\label{INV}Este teorema consta de 3 numerales:\newline
1. Si $A$ tiene la descomposici\'{o}n can\'{o}nica de Jordan%
\begin{equation}
A=T\left(
\begin{array}{cc}
C & 0 \\
0 & N%
\end{array}%
\right) T^{-1}.  \label{R9.0}
\end{equation}%
donde $C$ $\in
\mathbb{C}
^{p\times p}$ y $T\in
\mathbb{C}
^{m\times m}$ son matrices invertibles, $N\in
\mathbb{C}
^{q\times q}$ es una matriz nilpotente con \'{\i}ndice $k$ con $\left(
p,q\right) $ en $%
\mathbb{N}
^{2}\ $cumpliendo $p+q=n,$ entonces la inversa Drazin de $A$, $A^{D}$, se
obtiene de la siguiente forma%
\begin{equation}
A^{D}=T\left(
\begin{array}{cc}
C^{-1} & 0 \\
0 & 0%
\end{array}%
\right) T^{-1}  \label{R9}
\end{equation}%
2. Para cierto polinomio $Q\left( x\right) ,$ La inversa Drazin de $A$ puede
ser escrita como
\begin{equation*}
A^{D}=Q\left( A\right)
\end{equation*}%
3. Si en la descomposici\'{o}n can\'{o}nica de $A$ en (\ref{R9.0}) $N=0$,
entonces%
\begin{equation*}
A^{G}=T\left(
\begin{array}{cc}
C^{-1} & 0 \\
0 & 0%
\end{array}%
\right) T^{-1}.
\end{equation*}%
Demostraci\'{o}n: V\'{e}ase \cite{P1} y \cite{P340a}.
\end{theorem}

\begin{definition}
Sea $A\in
\mathbb{C}
^{m\times m},$ $f:%
\mathbb{C}
\rightarrow
\mathbb{C}
,$ decimos que $f$ pertece al espacio $\Im \left( A\right) $ si existe un
entorno $V$ del espectro de $A,$ $\sigma \left( A\right) $ sobre el cual f
es anal\'{\i}tica.

\begin{theorem}
\label{t0}$\left( Aplicaci\acute{o}n\text{ }espectral\right) $ Sea $f$ $\in
\Im \left( A\right) $, entonces
\begin{equation*}
f\left( \sigma \left( A\right) \right) =\sigma \left( f\left( A\right)
\right) ,
\end{equation*}%
donde $f\left( \sigma \left( A\right) \right) =\left\{ f\left( \lambda
\right) :\lambda \in \sigma \left( A\right) \right\} .$\newline
Demostraci\'{o}n: V\'{e}ase \cite{Dunford}.
\end{theorem}
\end{definition}

\begin{definition}
Sea $x\in
\mathbb{C}
^{m}.$ definimos la norma vectorial $\left\Vert {}\right\Vert _{1}$ por: $%
\left\Vert x\right\Vert _{1}=\left\Vert \left( x_{1},...,x_{m}\right)
\right\Vert _{1}=\sum\limits_{i=1}^{m}\left\vert x_{i}\right\vert ,$ y la
norma matricial por ella inducida como%
\begin{equation}
\left\Vert A\right\Vert _{1}=\sup \left\{ \left\Vert Ax\right\Vert _{1}:x\in
\mathbb{C}
^{m},\text{ }\left\Vert x\right\Vert _{1}=1\right\} ,\text{ \ para }A\in
\mathbb{C}
^{m\times m}  \label{xz8}
\end{equation}%
si las columnas de $A$ son $a_{1},$ $a_{2},...,a_{n}$ para las cuales
definimos sus componentes $a_{ij},$ $i=1,...,n$, para $j=1,...,n$
respectivamente, entonces $\left( \text{\ref{xz8}}\right) $ cohincide con
\begin{equation*}
\left\Vert A\right\Vert _{1}=\max \left\{ \sum\limits_{i=1}^{m}\left\vert
a_{ij}\right\vert :j=1,...,n\right\}
\end{equation*}
\end{definition}

\begin{theorem}
\label{AE1}Sea $A\in
\mathbb{C}
^{m\times m},$ $f$ $\in \Im \left( A\right) ,$ la funci\'{o}n matricial f(A)
se define por la matriz.%
\begin{equation*}
f(A)=P(A)
\end{equation*}%
donde P(x) es un polinomio conocido de grado menor que el del polinomio
minimal de A.\newline
Demostraci\'{o}n: V\'{e}ase \cite{Dunford}.
\end{theorem}

\begin{theorem}
\label{sig}Sea $A\in
\mathbb{C}
^{m\times m}$ talque $a=max\left\{ \left\vert z\right\vert ;z\in \sigma
\left( A\right) \right\} $, para $\epsilon >0$ existe una norma matricial en
$%
\mathbb{C}
^{m\times m}$ talque $\left\Vert A\right\Vert \leq a+\epsilon $.\newline
Dada la variedad de tipos de estabilidad manejadas en la literatura de los
esquemas de diferencias finitos num\'{e}ricos asociados a ecuaciones en
derivadas parciales, nosotros hemos optado por la siguiente definici\'{o}n:
\end{theorem}

\begin{definition}
\cite{P10a}$\left( Estabilidad\right) $. Consideremos fijo el par ordenado $%
\left( h,T\right) \in
\mathbb{R}
^{+}\times
\mathbb{R}
^{+}.$ Si se tiene una malla de puntos en el plano $\left( x,t\right) $ con
paso espacial $h$ y paso temporal $k>0$ talque $M$ es un n\'{u}mero natural.
Si $U$ es una funci\'{o}n definida en dicha malla y $U\left( ih,jk\right) $
representa su valor en el nodo $\left( i,j\right) ,$ se dice que $U$ es
estable, cuando $U\left( ih,jk\right) $ permanece acotada independiente de
los valores de $k$ y $M$ tal que $kM=T$
\end{definition}

\begin{definition}
(Problema de Sturm Liouville Discreto) El problema de valor en la frontera
formado por la ecuaci\'{o}n en diferencias
\begin{equation}
\Delta \left( p\left( i-1\right) \right) \Delta u\left( i-1\right) +q\left(
i\right) u\left( i\right) +\lambda r\left( i\right) u\left( i\right) =0,%
\text{ \ }0<i<N  \label{Q1}
\end{equation}%
con condiciones de frontera%
\begin{equation}
u\left( 0\right) =\alpha u\left( 1\right) ,\text{ \ \ \ \ \ \ \ \ \ \ \ }%
u\left( N\right) =\beta u\left( N-1\right) ,  \label{Q2}
\end{equation}%
en las que las funciones son escalares se llama problema de Sturm-Liouville
discreto (PSD).\newline
Donde $\Delta $ es el operador de diferencia definido por%
\begin{equation*}
\Delta p\left( i-1\right) =p\left( i\right) -p\left( i-1\right) ,\text{ \ \
\ }\Delta ^{2}p\left( i\right) =\Delta \left( \Delta p\left( i\right)
\right) .
\end{equation*}%
en la ecuaci\'{o}n en diferencias $\left( \text{\ref{Q1}}\right) $, $\lambda
$ es un par\'{a}metro, $\alpha $ y $\beta $ son constantes conocidas, la
funci\'{o}n $p\left( i\right) $ esta de finida en $0\leq i\leq N-1,$ $%
p\left( i\right) ,$ $q\left( i\right) $ y $r\left( i\right) $ en $0<i<N;$
adem\'{a}s $p\left( i\right) >0$ y $r\left( i\right) >0$.\newline
Al efectuar la operaci\'{o}n $\Delta $ en la ecuaci\'{o}n en diferencia $%
\left( \text{\ref{Q1}}\right) $, encontramos para $0<i<N$%
\begin{equation*}
p\left( i\right) u\left( i+1\right) -\left( p\left( i\right) +p\left(
i-1\right) \right) u\left( i\right) +\left( q\left( i\right) +\lambda
r\left( i\right) \right) u\left( i\right) +p\left( i-1\right) u\left(
i-1\right) =0
\end{equation*}%
denotando $s\left( i\right) =p\left( i\right) +p\left( i-1\right) -q\left(
i\right) ,$ $0<i<N$, se tiene%
\begin{equation}
-p\left( i-1\right) u\left( i-1\right) +s\left( i\right) u\left( i\right)
-p\left( i\right) u\left( i+1\right) =\lambda r\left( i\right) u\left(
i\right) ,\text{ }0<i<N  \label{R5}
\end{equation}%
en esta \'{u}ltima ecuaci\'{o}n para los dos casos\ $k=1$, $N-1$ tenemos las
ecuaciones%
\begin{eqnarray*}
-p\left( 0\right) u\left( 0\right) +s\left( 1\right) u\left( 1\right)
-p\left( 1\right) u\left( 2\right) &=&\lambda r\left( 1\right) u\left(
1\right) ,\text{ } \\
-p\left( N-2\right) u\left( N-2\right) +s\left( N-1\right) u\left(
N-1\right) -p\left( N-1\right) u\left( N\right) &=&\lambda r\left(
N-1\right) u\left( N-1\right) ,
\end{eqnarray*}%
las cuales teniendo en cuenta $\left( \text{\ref{Q2}}\right) $ se reescriben
como%
\begin{eqnarray}
\bar{s}\left( 1\right) u\left( 1\right) -p\left( 1\right) u\left( 2\right)
&=&\lambda r\left( 1\right) u\left( 1\right) ,  \label{R3} \\
-p\left( N-2\right) u\left( N-2\right) +\bar{s}\left( N-1\right) u\left(
N-1\right) &=&\lambda r\left( N-1\right) u\left( N-1\right) ,\text{ }
\label{R4}
\end{eqnarray}%
donde $\bar{s}\left( 1\right) =s\left( 1\right) -\alpha p\left( 0\right) $ y
$\bar{s}\left( N-1\right) =s\left( N-1\right) -\beta p\left( N-1\right) .$%
\newline
Con las ecuaciones $\left( \text{\ref{R3}}\right) $ para $i=0,$ $\left(
\text{\ref{R4}}\right) $ para $i=N-1$ y $\left( \text{\ref{R5}}\right) $
para $2\leq i\leq N-2$ se producen $N-1$ ecuaciones equivalentes a $\left(
\text{\ref{Q1}}\right) -\left( \text{\ref{Q2}}\right) $, las cuales
describen un problema generalizado de valores propios matricial con la forma
siguiente%
\begin{equation}
\begin{array}{c}
\left[
\begin{array}{cccccc}
\bar{s}\left( 1\right) & -p\left( 1\right) & 0 &  &  & 0 \\
-p\left( 1\right) & s\left( 2\right) & -p\left( 2\right) &  &  & \vdots \\
0 & -p\left( 2\right) & \ddots &  &  &  \\
\vdots &  &  &  & s\left( N-3\right) & -p\left( N-2\right) \\
0 &  &  &  & -p\left( N-2\right) & \bar{s}\left( N-1\right)%
\end{array}%
\right] \left[
\begin{array}{c}
u\left( 1\right) \\
u\left( 2\right) \\
\vdots \\
u\left( N-2\right) \\
u\left( N-1\right)%
\end{array}%
\right] \\
=\lambda \left[
\begin{array}{cccccc}
r\left( 1\right) & 0 & 0 &  &  & 0 \\
0 & r\left( 2\right) &  &  &  & \vdots \\
\vdots & 0 & \ddots &  &  & 0 \\
0 &  &  &  & r\left( N-2\right) & 0 \\
0 &  &  &  & 0 & r\left( N-1\right)%
\end{array}%
\right] \left[
\begin{array}{c}
u\left( 1\right) \\
u\left( 2\right) \\
\vdots \\
u\left( N-2\right) \\
u\left( N-1\right)%
\end{array}%
\right]%
\end{array}
\label{ss33}
\end{equation}%
Ya que $p\left( i\right) >0$, $0\leq i\leq N-1;$ $r\left( i\right) >0$, $%
0<i<N$ de acuerdo a la teoria espectral para matrices sim\'{e}tricas
nosotros asumimos los siguientes teoremas.
\end{definition}

\begin{theorem}
\label{TT1}El PSD tiene exactamente $N-1$ eigenvalores $\lambda _{l},$ $%
0<l<N $; a cada uno de los cuales les corresponde las eigenfunciones $%
v_{l}(i),$ $0<i<N$ definidas para $0\leq i\leq N$.\newline
Demostraci\'{o}n: V\'{e}ase \cite{Argawal} \'{o} \cite{P10a}.
\end{theorem}

\begin{theorem}
Sean $\lambda _{l}$ y $v_{l}(i)$ los $N-1$ eigenvalores y egenfuciones del
PSD, entonces el conjunto
\begin{equation*}
\left\{ v_{l}(i):l=1,...,N-1\right\}
\end{equation*}%
es ortogonal con respecto al producto interno
\begin{equation*}
\sum_{i=1}^{N-1}r\left( i\right) v_{\mu }\left( i\right) v_{\upsilon }\left(
i\right) =0\text{, \ \ \ \ \ \ \ }\forall \mu \neq \upsilon .
\end{equation*}%
Demostraci\'{o}n: V\'{e}ase \cite{Argawal}.
\end{theorem}

\begin{theorem}
\label{TTT1}Sea $\left\{ v_{l}(i):0<l<N\right\} $ el conjunto de
eigenfunciones en el teorema anterior y $u\left( i\right) $ cualquier funci%
\'{o}n definida para $0<i<N\ $, entonces $u\left( i\right) $ puede
expresarse como una combinaci\'{o}n lineal de dicho conjunto$,$ es decir:%
\begin{equation*}
u\left( i\right) =\sum_{l=1}^{N-1}c_{l}v_{l}\left( i\right) ,\text{\ \ \ \ }%
0<i<N
\end{equation*}%
donde cada $c_{l}$ puede ser determinada por%
\begin{equation*}
c_{l}=\frac{\sum_{i=1}^{N-1}r\left( i\right) v_{l}\left( i\right) u\left(
i\right) }{\sum_{i=1}^{N-1}r\left( i\right) v_{l}^{2}\left( i\right) },\text{%
\ \ \ \ }0<l<N
\end{equation*}%
Demostraci\'{o}n: V\'{e}ase \cite{Argawal}.\newline
Las siguientes 3 proposiciones las enunciamos con una notaci\'{o}n
pertinente para la demostraci\'{o}n de nuestro teorema fundamental $\left(
\text{\ref{Main}}\right) $.
\end{theorem}

\begin{proposition}
$:$\label{T1 copy(2)}Sean $\alpha $ y $\beta $ escalares conocidos, $N\in
\mathbb{N}
$ y $\lambda $ un eigenvalor del PSD%
\begin{eqnarray}
h\left( i+1\right) -\left( 2-\lambda \right) h\left( i\right) +h\left(
i-1\right) &=&0,\text{ }0<i<N,  \label{32a} \\
h\left( 0\right) +\alpha N\left( h\left( 1\right) -h\left( 0\right) \right)
&=&0,  \label{33a} \\
h\left( N\right) +\beta N\left( h\left( N\right) -h\left( N-1\right) \right)
&=&0,  \label{34a}
\end{eqnarray}%
con eigenfunci\'{o}n $v\left( i\right) .$ Entonces $v\left( i\right) $ es
soluci\'{o}n del PSD $\left( \ref{32a}\right) -\left( \ref{34a}\right) $ si
y s\'{o}lo si%
\begin{equation*}
H(i)=v\left( i\right) R
\end{equation*}%
con $R$ vector arbitrario en $%
\mathbb{C}
^{m},$ es una soluci\'{o}n del problema de Sturm Liouville discreto
vectorial (PSDV)%
\begin{eqnarray}
H\left( i+1\right) -\left( 2-\lambda \right) H\left( i\right) +H\left(
i-1\right) &=&0,\text{ }0<i<N,  \label{32b} \\
H\left( 0\right) +\alpha N\left( H\left( 1\right) -H\left( 0\right) \right)
&=&0,  \label{33b} \\
H\left( N\right) +\beta N\left( H\left( N\right) -H\left( N-1\right) \right)
&=&0,  \label{34b}
\end{eqnarray}%
en el que $H\left( i\right) $ es una funci\'{o}n vectorial discreta en $%
\mathbb{C}
^{m}.$\newline
Demostraci\'{o}n:\linebreak Supongamos que $H\left( i\right) $ satisface el
PSDV $\left( \ref{32b}\right) -\left( \ref{34b}\right) $ donde $H\left(
i\right) =\left(
\begin{array}{cccc}
h_{1}\left( i\right) & h_{2}\left( i\right) & ... & h_{m}\left( i\right)%
\end{array}%
\right) $ siendo $h_{k}\left( i\right) $ funciones escalares discretas
definidas para $0<i<N.$ Luego en virtud de la igualdad entre vectores
tenemos que cada una de sus componentes $h_{k}\left( i\right) $ $k=1,...,m$
satisfacen el PSD $\left( \ref{32a}\right) -\left( \ref{34a}\right) .$
Entonces al hallar una soluci\'{o}n para una componente, ella lo ser\'{a}
tambi\'{e}n para las dem\'{a}s. Sin p\'{e}rdida de generalidad denotemos $%
h_{k}\left( i\right) =h\left( i\right) $. Si en $\left( \ref{ss33}\right) $
definimos:
\begin{equation*}
q(i)=0,\text{ }r(i)=1,\text{ }p(i)=1,
\end{equation*}%
de $\left( \ref{32a}\right) -\left( \ref{34a}\right) $ se obtiene
\begin{equation}
\left[
\begin{array}{cccccc}
\frac{2-\alpha N}{1-\alpha N} & -1 & 0 &  &  & 0 \\
-1 & 2 & -1 &  &  & \vdots \\
0 & -1 & \ddots &  &  &  \\
\vdots &  &  &  & 2 & -1 \\
0 &  &  &  & -1 & \frac{2+\beta N}{1+\beta N}%
\end{array}%
\right] \left[
\begin{array}{c}
h\left( 1\right) \\
h\left( 2\right) \\
\vdots \\
h\left( N-2\right) \\
h\left( N-1\right)%
\end{array}%
\right] =\lambda \left[
\begin{array}{c}
h\left( 1\right) \\
h\left( 2\right) \\
\vdots \\
h\left( N-2\right) \\
h\left( N-1\right)%
\end{array}%
\right]  \label{STL}
\end{equation}%
Por el teorema $\left( \ref{TT1}\right) ,$ $\left( \ref{STL}\right) $ tiene
los $N-1$ eigenvalores%
\begin{equation}
\left\{ \lambda _{l}\right\} _{l=1}^{N-1}  \label{eige1}
\end{equation}%
con eigenfunciones asociadas
\begin{eqnarray}
&&\text{ \ \ \ \ \ \ \ \ \ }%
\begin{array}{c}
\left\{ v_{l}(i)\right\} _{l=1}^{N-1},\text{ \ }0<l<N,\text{ }0<i<N,%
\end{array}
\label{eigenf} \\
&&%
\begin{array}{c}
v_{l}=\left( v_{l}\left( 1\right) ,\text{ }v_{l}\left( 2\right) ,\text{..}%
v_{l}\left( N-1\right) \right) ^{T},\text{ \ }0<l<N.%
\end{array}
\notag
\end{eqnarray}%
Luego, las soluciones del PSD $\left( \ref{32a}\right) -\left( \ref{34a}%
\right) $ tienen la forma
\begin{equation*}
v_{l}(i)=v_{l}\left( i\right) D_{2},\text{ con }D_{2}\text{ arbitrario en }%
\mathbb{C}\text{ e }i=1,...N-1.
\end{equation*}%
Como lo anterior se verifica para cada componente de $H(i)$ concluimos que%
\begin{equation}
H(i)=(h_{1}\left( i\right) ,...,h_{m}\left( i\right) )=v_{l}\left( i\right)
R,\text{ \ }R\text{ arbitrario en }\mathbb{C}^{m}  \label{TT2}
\end{equation}%
es soluci\'{o}n del PSDV $\left( \ref{32b}\right) -\left( \ref{34b}\right) $.%
\newline
Aplicando el teorema $\left( \text{\ref{TT1}}\right) $ es claro afirmar que
el PSDV $\left( \ref{32b}\right) -\left( \ref{34b}\right) $ tiene
exactamente $N-1$ soluciones de la forma $\left( \text{\ref{TT2}}\right) $
donde $v_{l}\left( i\right) $ son las $N-1$ eigenfunciones de PSD $\left( %
\ref{32a}\right) -\left( \ref{34a}\right) $ determinadas en $\left( \ref%
{eigenf}\right) .$
\end{proposition}

\begin{proposition}
$:$\label{SL} Sea $H\left( i\right) $ funci\'{o}n vectorial discreta en $%
\mathbb{C}
^{m}$, que verifica%
\begin{eqnarray}
H\left( i+1\right) -\left( 2-\lambda \right) H\left( i\right) +H\left(
i-1\right) &=&0,\text{ }0<i<N,  \label{D1} \\
A_{1}H\left( 0\right) +A_{2}N\left( H\left( 1\right) -H\left( 0\right)
\right) &=&0,  \label{D2} \\
B_{1}H\left( N\right) +B_{2}N\left( H\left( N\right) -H\left( N-1\right)
\right) &=&0,  \label{D3}
\end{eqnarray}%
con $A_{1}$, $A_{2},$ $B_{1},$ $B_{2}\in \mathbb{C}$ $^{m\times m}$ y%
\begin{equation}
G\left( \alpha ,\beta \right) =\left(
\begin{array}{c}
\alpha A_{1}-A_{2} \\
\beta B_{1}-B_{2}%
\end{array}%
\right) ,  \label{RR}
\end{equation}%
tal que $rang\left( \tilde{G}\left( \alpha ,\beta \right) \right) <m$ y $%
R\in KerG\left( \alpha ,\beta \right) $. Entonces para cada $\lambda _{l},$ $%
0<l<N$ eigenvalor del PSD $\left( \ref{32a}\right) -\left( \ref{34a}\right) $%
\begin{equation*}
H_{l}(i)=v_{l}\left( i\right) R,\text{ }0<i<N
\end{equation*}%
son $N-1$ soluciones de $\left( \ref{D1}\right) -\left( \ref{D2}\right) .$%
\newline
Demostraci\'{o}n: Sea
\begin{equation*}
H_{l}(i)=v_{l}\left( i\right) R\text{ }0<i<N\text{,\ }R\in KerG\left( \mu
,\beta \right)
\end{equation*}%
con $H(i)$ as\'{\i} definida, $\left( \ref{D1}\right) -\left( \ref{D2}%
\right) $ es equivalente a
\begin{eqnarray}
\left( v\left( i+1\right) -\left( 2-\lambda \right) v\left( i\right)
+v\left( i-1\right) \right) R &=&0,\text{ }0<i<N,  \label{31ab} \\
\left( A_{1}v\left( 0\right) +A_{2}N\left( v\left( 1\right) -v\left(
0\right) \right) \right) R &=&0,  \label{32ab} \\
\left( B_{1}v\left( N\right) +B_{2}N\left( v\left( N\right) -v\left(
N-1\right) \right) \right) R &=&0,  \label{33ab}
\end{eqnarray}%
como $v\left( i\right) $ verifica $\left( \ref{32a}\right) -\left( \ref{34a}%
\right) ,$ se tiene que $\left( \ref{31ab}\right) -\left( \ref{33ab}\right)
, $ es equivalente a%
\begin{eqnarray*}
\left( v\left( i+1\right) -\left( 2-\lambda \right) v\left( i\right)
+v\left( i-1\right) \right) R &=&0,\text{ }0<i<N, \\
\left( \alpha A_{1}-A_{2}\right) (v\left( 1\right) -v\left( 0\right) )R &=&0,
\\
\left( \beta B_{1}-B_{2}\right) (v\left( N\right) -v\left( N-1\right) )R
&=&0.
\end{eqnarray*}%
Como adem\'{a}s $R\in KerG\left( \alpha ,\beta \right) $, se concluye que $%
H_{l}(i)=v_{l}\left( i\right) R$, satisface $\left( \ref{D1}\right) -\left( %
\ref{D2}\right) $.\newline
El teorema y proposici\'{o}n siguientes los enunciamos con una notaci\'{o}n
pertinente para la demostraci\'{o}n de nuestro teorema fundamental.
\end{proposition}

\begin{theorem}
\label{AA} Sean $I,$ $A\in
\mathbb{C}
^{m\times m}$ donde $I$ es la identidad, y $G\left( j\right) $ funci\'{o}n
vectorial definida para $i=0,...N;$ supongamos adem\'{a}s $\rho $ escalar y
consideremos la siguiente ecuaci\'{o}n en diferencias matricial
\begin{equation*}
G\left( j+1\right) -(2I+\rho A)G\left( j\right) +G\left( j-1\right) =0,\text{%
\ }j>0.\text{ }G\left( j\right) \in \mathbb{C}^{m\times m}
\end{equation*}%
Entonces su soluci\'{o}n general es%
\begin{equation*}
G\left( j\right) =\left[ P_{+}\left( A\right) \right] ^{j}P_{1}+\left[
P_{-}\left( A\right) \right] ^{j}Q_{1},\text{\ }j>0,\text{ }l_{1},l_{2}\text{
}\in \mathbb{C}^{m},
\end{equation*}%
donde $P_{+}\left( x\right) $ y $P_{-}\left( x\right) $ son los polinomios
de grado m\'{a}s peque\~{n}o que el grado del polinomio minimal de $A$ tal
que
\begin{equation*}
P_{+}\left( A\right) =I+\frac{\rho }{2}A+\sqrt{\left( I+\frac{\rho }{2}%
A\right) ^{2}-I},P_{-}\left( A\right) =I+\frac{\rho }{2}A-\sqrt{\left( I+%
\frac{\rho }{2}A\right) ^{2}-I}.
\end{equation*}%
Demostraci\'{o}n: V\'{e}ase \cite{camacho}, o tambi\'{e}n \cite{P37a}.
\end{theorem}

\begin{proposition}
\label{AAA} Con las mismas hip\'{o}tesis del teorema anterior, consideremos
la siguiente ecuaci\'{o}n en diferencias matricial%
\begin{equation}
\begin{array}{c}
EG\left( j+1\right) -(2E+\rho A)G\left( j\right) +EG\left( j-1\right) =0,%
\text{\ }j>0.\text{ }%
\end{array}
\label{16a}
\end{equation}%
donde $E\in
\mathbb{C}
^{m\times m}$\ es singular, supongamos adem\'{a}s que existe $\gamma \in
\mathbb{C}
$ tal que $(\gamma E+A)$ es no singular, defininamos $\hat{E}$ $=(\gamma
E+A)^{-1}E$, $\hat{A}=(\lambda E+A)^{-1}A$. Entonces su soluci\'{o}n general
viene dada por%
\begin{eqnarray*}
G\left( j\right) &=&\hat{Z}_{0}^{\text{ \ }j}\text{\ }\hat{E}\hat{E}%
^{D}l_{1}+\hat{Z}_{1}^{\text{ \ }j}\hat{E}\hat{E}^{D}l_{2}\text{, \ con \ }%
l_{1},\text{ }l_{2}\text{ arbitrarios en }%
\mathbb{C}
^{m}. \\
\text{donde }\hat{Z}_{0}^{\text{ }} &=&\left[ P_{+}\left( \hat{E}^{D}\hat{A}%
\right) \right] \hat{E}\hat{E}^{D},\text{ }\hat{Z}_{1}^{\text{ }}=\left[
P_{-}\left( \hat{E}^{D}\hat{A}\right) \right] \hat{E}\hat{E}^{D}
\end{eqnarray*}%
con \
\begin{equation*}
P_{+}\left( \hat{E}^{D}\hat{A}\right) =I+\frac{\rho }{2}\hat{E}^{D}\hat{A}+%
\sqrt{\left( I+\frac{\rho }{2}\hat{E}^{D}\hat{A}\right) ^{2}-I},\text{ }%
P_{-}\left( \hat{E}^{D}\hat{A}\right) =I+\frac{\rho }{2}\hat{E}^{D}\hat{A}-%
\sqrt{\left( I+\frac{\rho }{2}\hat{E}^{D}\hat{A}\right) ^{2}-I}
\end{equation*}%
Demostraci\'{o}n: Dado que existe $\gamma $ con $(\gamma E+A)$ invertible,
sean $\hat{E}$ $=(\gamma E+A)^{-1}E$, $\hat{A}=(\lambda E+A)^{-1}A,$ luego
la ecuaci\'{o}n $\left( \ref{16a}\right) $ es equivalente a%
\begin{equation}
\hat{E}G\left( j+1\right) -(2\hat{E}-\rho \hat{A})G\left( j\right) +\hat{E}%
G\left( j-1\right) =0,\text{\ }j>0.  \label{23}
\end{equation}%
Considerando la descomposici\'{o}n can\'{o}nica de $\hat{E}$, obtenemos:
\begin{equation}
\hat{E}=T^{-1}\left(
\begin{array}{cc}
C & 0 \\
0 & N%
\end{array}%
\right) T,\text{ }\hat{A}=T^{-1}\left(
\begin{array}{cc}
I-\gamma C & 0 \\
0 & I-\gamma N%
\end{array}%
\right) T  \label{rr22}
\end{equation}%
en las que $C$ $\in
\mathbb{C}
^{p\times p}$ es matriz invertible, $N\in
\mathbb{C}
^{q\times q}$ es una matriz nilpotente de \'{\i}ndice $k$ con $p$, $q\ $%
cumpliendo $p+q=n$. Luego haciendo
\begin{equation}
TG\left( j\right) =\left(
\begin{array}{c}
h\left( j\right) \\
d\left( j\right)%
\end{array}%
\right) ,  \label{AEE1}
\end{equation}%
$\left( \ref{16a}\right) $ es equivalente a los dos siguientes sistemas:%
\begin{equation}
Ch\left( j+1\right) -\left[ 2C+\rho \left( I-\gamma C\right) \right] h\left(
j\right) +Ch\left( j-1\right) =0,  \label{24}
\end{equation}%
\begin{equation}
Nd\left( j+1\right) -\left[ 2N+\rho \left( I-\gamma N\right) \right] d\left(
j\right) +Nd\left( j-1\right) =0.  \label{25}
\end{equation}%
Como $C$ es invertible se obtiene que $\left( \ref{24}\right) $ se verifica
si y s\'{o}lo si
\begin{equation*}
h\left( j+1\right) -\left[ 2I+\rho C^{-1}\left( I-\gamma C\right) \right]
h\left( j\right) +h\left( j-1\right) =0,
\end{equation*}%
la cual por el teorema \ref{AA}, tiene la soluci\'{o}n general
\begin{equation*}
h\left( j\right) =\left[ P_{+}\left( D\right) \right] ^{j}P_{1}+\left[
P_{-}\left( D\right) \right] ^{j}Q_{1},\text{\ }j>0,\text{ }l_{1},l_{2}\text{
}\in \mathbb{C}^{k}.
\end{equation*}%
donde $D=C^{-1}\left( I-\gamma C\right) $ y $P_{+}\left( x\right) $ y $%
P_{-}\left( x\right) $ son los polinomios de grado m\'{a}s peque\~{n}o que
el grado del polinomio minimal de $D$ tal que
\begin{equation*}
P_{+}\left( D\right) =I+\frac{\rho }{2}D+\sqrt{\left( I+\frac{\rho }{2}%
D\right) ^{2}-I},\text{ }P_{-}\left( D\right) =I+\frac{\rho }{2}D-\sqrt{%
\left( I+\frac{\rho }{2}D\right) ^{2}-I}
\end{equation*}%
Como para N suficiente mente grande $2N+\rho \left( I-\gamma N\right) $ es
invertible y $k$ es el indice de nilpotencia de $N$, por inducci\'{o}n puede
demostrarse que la ecuaci\'{o}n $\left( \ref{25}\right) $ tiene la soluci%
\'{o}n trivial $d\left( j\right) =0$, $j>0$, luego considerando
\begin{equation}
\hat{E}=T^{-1}\left(
\begin{array}{cc}
C & 0 \\
0 & N%
\end{array}%
\right) T,\text{ }\hat{E}^{D}=T^{-1}\left(
\begin{array}{cc}
C^{-1} & 0 \\
0 & 0%
\end{array}%
\right) T  \label{ez2}
\end{equation}%
y haciendo uso de las ecuaciones $\left( \ref{AEE1}\right) ,$ $\left( \ref%
{24}\right) $ y $\left( \ref{25}\right) ,$ se concluye que $\left( \ref{23}%
\right) $ tiene la soluci\'{o}n general%
\begin{eqnarray}
G\left( j\right) &=&\hat{Z}_{0}^{\text{ \ }j}\text{\ }\hat{E}\hat{E}%
^{D}l_{1}+\hat{Z}_{1}^{\text{ \ }j}\hat{E}\hat{E}^{D}l_{2}\text{, \ con \ }%
l_{1},\text{ }l_{2}\text{ arbitrarios en }%
\mathbb{C}
^{m}.  \label{r.2} \\
\text{donde }\hat{Z}_{0}^{\text{ }} &=&\left[ P_{+}\left( \hat{E}^{D}\hat{A}%
\right) \right] \hat{E}\hat{E}^{D},\text{ }\hat{Z}_{1}^{\text{ }}=\left[
P_{-}\left( \hat{E}^{D}\hat{A}\right) \right] \hat{E}\hat{E}^{D}
\label{r.2b}
\end{eqnarray}
\end{proposition}

\section{Soluci\'{o}n discretizada de la ecuaci\'{o}n de onda}

Retomemos nuestro problema central

\begin{eqnarray}
Eu_{tt}(x,t)-Au_{xx}(x,t) &=&0,\ x\in \lbrack 0,1],\text{ }t>0,  \label{1} \\
A_{1}u(0,t)+A_{2}u_{x}(0,t) &=&0,\ t>0,  \label{2} \\
B_{1}u(1,t)+B_{2}u_{x}(1,t) &=&0,\text{ }t>0,  \label{2.1} \\
u(x,0) &=&f(x),\ x\in \lbrack 0,1],\text{ \ }  \label{3} \\
u_{t}(x,0) &=&g(x),\text{ }x\in \lbrack 0,1],  \label{4}
\end{eqnarray}%
En el que asumimos adem\'{a}s la hipotesis: que existe $\gamma \in
\mathbb{C}
$ tal que $(\gamma E+A)$ es no singular, y as\'{\i} definimos $\hat{E}$ $%
=(\gamma E+A)^{-1}E$ , $\hat{A}=(\lambda E+A)^{-1}A$ con la condici\'{o}n:%
\begin{equation}
z\neq 0\text{ para alg\'{u}n }z\in \sigma \left( \hat{E}^{D}\hat{A}\right)
\subset \mathbb{C}  \label{5}
\end{equation}%
En el esquema de diferencias correspondiente a $\left( \ref{1}\right)
-\left( \ref{4}\right) $ dividimos el dominio $\left[ 0,1\right] \times %
\left] 0,\infty \right[ $ en rect\'{a}ngulos de lados $\Delta x=h,$ $\Delta
t=k$, luego al introducir coordenadas de un punto t\'{\i}pico de la malla $%
\left( ih,jk\right) ,$ representamos el valor $u\left( ih,jk\right) $ por $%
U\left( i,j\right) .$ Si aproximamos las derivadas parciales de $u$
utilizando diferencias avanzadas para las derivadas primeras y diferencias
centrales para las derivadas segundas. Ellas tienen la forma%
\begin{eqnarray}
u_{t}\left( ih,jk\right) &=&\frac{U\left( i,j+1\right) -U\left( i,j\right) }{%
k},  \label{k1} \\
u_{x}\left( ih,jk\right) &=&\frac{U\left( i+1,j\right) -U\left( i,j\right) }{%
h},  \label{k2} \\
u_{xx}\left( ih,jk\right) &=&\frac{U\left( i+1,j\right) -2U\left( i,j\right)
+U\left( i-1,j\right) }{h^{2}},  \label{k3} \\
u_{tt}\left( ih,jk\right) &=&\frac{U\left( i,j+1\right) -2U\left( i,j\right)
+U\left( i,j-1\right) }{k^{2}}.  \label{k4}
\end{eqnarray}%
Al sustituir $\left( \ref{k1}\right) -\left( \ref{k4}\right) $ en $\left( %
\ref{1}\right) -\left( \ref{4}\right) ,$ para $N$ un entero positivo con $%
h=1/N$ y $r=k/h$ $0<i<N;$ se obtiene la siguiente representaci\'{o}n
aproximada de $\left( \ref{1}\right) -\left( \ref{4}\right) $

\begin{equation}
\begin{array}{c}
r^{2}A\left( U\left( i+1,j\right) -U\left( i-1,j\right) \right) +2\left(
I-r^{2}A\right) U\left( i,j\right) - \\
\left( U\left( i,j+1\right) +U\left( i,j-1\right) \right) =0,\text{ }0<i<N,%
\text{ }j>0%
\end{array}
\label{6}
\end{equation}%
\begin{gather}
A_{1}U(0,j)+NA_{2}\left( U(1,j)-U(0,j)\right) =0,\hspace{0.5cm}j>0,
\label{7} \\
B_{1}U(N,j)+NB_{2}\left( U(N,j)-U(N-1,j)\right) =0,\hspace{0.5cm}j>0,
\label{8} \\
U(i,0)=F(i),\text{ }0\leq i\leq N,  \label{9} \\
\dfrac{U(i,1)-U(i,0)}{k}=G(i),\text{ }0\leq i\leq n.  \label{10}
\end{gather}%
Busquemos soluciones no triviales a la ecuaci\'{o}n de onda discretizada $%
\left( \ref{6}\right) .$ Para esto usamos el m\'{e}todo de separaci\'{o}n de
variables discreto matricial suponiendo $U\left( i,j\right) $ de la forma%
\begin{equation}
U\left( i,j\right) =G\left( j\right) H\left( i\right) ,\text{ }H\left(
i\right) \in
\mathbb{C}
^{m}.\text{\ }  \label{12}
\end{equation}%
Para $r=k/h,$ imponiendo que $\left\{ U\left( i,j\right) \right\} $
verifique $\left( \ref{6}\right) $, resulta%
\begin{equation}
\begin{array}{c}
r^{2}AG\left( j\right) \left[ H\left( i+1\right) -H\left( i-1,j\right) %
\right] +2\left( E-Ar^{2}\right) G\left( j\right) H\left( i\right) - \\
E\left[ G\left( j+1\right) +U\left( j-1\right) H\left( i\right) \right] =0,%
\text{\ }0<i<N,\text{ }j>0%
\end{array}%
\text{ }  \label{13}
\end{equation}%
Tomando $\rho \in \mathbb{R}$ arbitrario y puesto que $2\left(
E-Ar^{2}\right) =-r^{2}A\left( 2+\frac{\rho }{r^{2}}\right) +\left( 2E+\rho
A\right) ,$ $\left( \ref{13}\right) $ toma la forma%
\begin{equation}
\begin{array}{c}
r^{2}AG\left( j\right) \left[ H\left( i+1\right) -A\left( 2+\frac{\rho }{%
r^{2}}\right) H\left( i\right) +H\left( i-1\right) \right] - \\
\left[ EG\left( j+1\right) -(2E+\rho AG\left( j\right) +EG\left( j-1\right) %
\right] H\left( i\right) =0,\text{ }0<i<N,\text{ }j>0.%
\end{array}
\label{14}
\end{equation}%
La ecuaci\'{o}n $\left( \ref{14}\right) $ se verifica si $\left\{ H\left(
i\right) \right\} $ y $\left\{ G\left( j\right) \right\} $ satisfacen las
ecuaciones en diferencias%
\begin{equation}
\begin{array}{c}
H\left( i+1\right) -\left( 2+\frac{\rho }{r^{2}}\right) H\left( i\right)
+H\left( i-1\right) =0,\text{ }0<i<N,%
\end{array}
\label{15}
\end{equation}%
\begin{equation}
\begin{array}{c}
EG\left( j+1\right) -(2E+\rho A)G\left( j\right) +EG\left( j-1\right) =0,%
\text{\ }j>0.\text{ }%
\end{array}
\label{16}
\end{equation}%
Si en $\left( \ref{15}\right) $ definimos $\rho =-r^{2}\lambda _{l}$ para
cada $\lambda _{l}$ eigenvalor del PSD $\left( \ref{32a}\right) -\left( \ref%
{34a}\right) ,$ de la proposici\'{o}n $\left( \text{\ref{T1 copy(2)}}\right)
,$ sabemos que si $v_{l}\left( i\right) $ es el eigenvector asociado a $%
\lambda _{l}.$ Entonces: $H_{l}\left( i\right) $ definida como

\begin{equation}
H_{l}\left( i\right) =v_{l}\left( i\right) R,\text{ }R\in
\mathbb{C}
^{m}-\left\{ 0\right\} \text{, \ }0<i<N,\text{ \ }j>0.  \label{rra}
\end{equation}%
satisface $\left( \ref{15}\right) $.\newline
Por otra parte usando la propocisi\'{o}n $\left( \ref{AAA}\right) $, para $%
\rho =-r^{2}\lambda _{l},$ la ecuaci\'{o}n $(\ref{16}),$ tiene soluci\'{o}n
dada por%
\begin{eqnarray}
G\left( j\right) &=&\hat{Z}_{0}^{\text{ \ }j}\text{\ }\hat{E}\hat{E}%
^{D}l_{1}+\hat{Z}_{1}^{\text{ \ }j}\hat{E}\hat{E}^{D}l_{2}\text{, \ con \ }%
l_{1},\text{ }l_{2}\text{ arbitrarios en }%
\mathbb{C}
^{m}.  \label{r.2} \\
\text{donde }\hat{Z}_{0}^{\text{ }} &=&\left[ P_{+}\left( \hat{E}^{D}\hat{A}%
\right) \right] \hat{E}\hat{E}^{D},\text{ }\hat{Z}_{1}^{\text{ }}=\left[
P_{-}\left( \hat{E}^{D}\hat{A}\right) \right] \hat{E}\hat{E}^{D}
\label{r.2b}
\end{eqnarray}%
as\'{\i}, una soluci\'{o}n de $\left( \ref{6}\right) $ es de la forma
\begin{equation}
U\left( i,j\right) =\hat{Z}_{0}^{\text{ \ }j}H_{1l}\left( i\right) +\hat{Z}%
_{0}^{\text{ \ }j}H_{2l}\left( i\right) ,\text{ \ \ }0<i<N,\text{ }j>0.
\label{12a}
\end{equation}%
donde
\begin{equation}
H_{1l}\left( i\right) =P_{l}v\left( i\right) ,\text{ }H_{2l}\left( i\right)
=Q_{l}v\left( i\right) \text{\ con\ }P_{l},\text{\ }Q_{l}\in
\mathbb{C}
^{m}-\left\{ 0\right\} .  \label{rra1}
\end{equation}

\subsection{Condici\'{o}n de contorno discretizada}

Las soluciones $U\left( i,j\right) $ de $\left( \ref{6}\right) $ encontradas
en $\left( \ref{12a}\right) $, se escriben de acuerdo a $\left( \ref{rra1}%
\right) ,$ como sigue%
\begin{equation}
U_{l}\left( i,j\right) =\hat{Z}_{0}^{\text{ \ }j}P_{l}v_{l}\left( i\right) +%
\hat{Z}_{1}^{\text{ \ }j}Q_{l}v_{l}\left( i\right) ,\text{ \ \ }0<i<N,\text{
}j>0.  \label{rar}
\end{equation}%
donde $v_{l}\left( i\right) $ es alguna eigenfunci\'{o}n del PSD $\left( \ref%
{32a}\right) -\left( \ref{34a}\right) $ y $P_{l},$\ $Q_{l}\in \left(
\mathbb{C}
^{m}-\left\{ 0\right\} \right) .$ Sea $G\left( \alpha ,\beta \right) $
definida en $\left( \ref{RR}\right) ,$ por la proposici\'{o}n $\left( \ref%
{SL}\right) ,$ dejando $j$ fijo, $\left( \ref{rar}\right) $ satisface $%
\left( \ref{7}\right) -\left( \ref{8}\right) $ si $P_{l}$ y $Q_{l}$ en los
sistemas de ecuaciones%
\begin{eqnarray*}
\hat{Z}_{0}^{\text{ \ }j}\left( A_{1}v\left( 0\right) +A_{2}N\left( v\left(
1\right) -v\left( 0\right) \right) \right) P_{l} &=&0, \\
\hat{Z}_{0}^{\text{ \ }j}\left( B_{1}v\left( N\right) +B_{2}N\left( v\left(
N\right) -v\left( N-1\right) \right) \right) P_{l} &=&0,
\end{eqnarray*}%
y%
\begin{eqnarray*}
\hat{Z}_{1}^{\text{ \ }j}\left( A_{1}v\left( 0\right) +A_{2}N\left( v\left(
1\right) -v\left( 0\right) \right) \right) Q_{l} &=&0, \\
\hat{Z}_{1}^{\text{ \ }j}\left( B_{1}v\left( N\right) +B_{2}N\left( v\left(
N\right) -v\left( N-1\right) \right) \right) Q_{l} &=&0,
\end{eqnarray*}%
pertenecen a $KerG\left( \mu ,\beta \right) \cap \left(
\mathbb{C}
^{m}-\left\{ 0\right\} \right) $. Con lo anterior $U_{l}\left( i,j\right) $
en $\left( \ref{rar}\right) $ es soluci\'{o}n de $\left( \ref{6}\right)
-\left( \ref{8}\right) .$\newline
Los resultados obtenidos hasta el momento, los resumimos en el siguiente
teorema.

\begin{theorem}
\label{Main0}Supongamos $f(x),$ $g(x)$ funciones en $%
\mathbb{C}
^{m},$ $A,$ $A_{1},$ $A_{2},$ $B_{1},$ $B_{2}$; $E$ $\in
\mathbb{C}
^{m\times m},$ $E$ matriz singular con $\gamma $ tal que $(\gamma E+A)$ es
invertible, $\hat{A},$ $\hat{E},$ $\hat{E}^{D}$ $\in
\mathbb{C}
$ $^{m\times m}$ definidas por $\left( \ref{rr22}\right) $, $\left( \ref{ez2}%
\right) $ y que verifican la condici\'{o}n $\left( \ref{5}\right) $; sean
adem\'{a}s $\left( \alpha ,\beta \right) \in
\mathbb{R}
^{2}$ que verifican $\left( \ref{32a}\right) -\left( \ref{34a}\right) ,$ $%
G\left( \alpha ,\beta \right) $ representada por $\left( \ref{RR}\right) $
que cumple $rang\left( G\left( \alpha ,\beta \right) \right) <m$. Entonces
una soluci\'{o}n de $\left( \ref{6}\right) -\left( \ref{8}\right) $ esta
dada por:%
\begin{equation}
U\left( i,j\right) =\sum_{l=1}^{N-1}\left( \hat{Z}_{0}^{\text{\ }j}\hat{E}%
\hat{E}^{D}P_{l}+\hat{Z}_{1}^{\text{\ }j}\hat{E}\hat{E}^{D}Q_{l}\right)
v_{l}\left( i\right) ,\text{ }0<i<N,\text{ }j>0  \label{ss}
\end{equation}%
donde $\hat{Z}_{0}$ y $\hat{Z}_{1}$ est\'{a}n definidas por $\left( \ref%
{r.2b}\right) $, $l_{1},$ $l_{2}\in KerG\left( \mu ,\beta \right) \cap
\left(
\mathbb{C}
^{m}-\left\{ 0\right\} \right) $ y$\ v_{l}\left( i\right) $ son las $N-1$
eigenfunciones asociadas al PSD $\left( \ref{32a}\right) -\left( \ref{34a}%
\right) $ determinadas seg\'{u}n la proposici\'{o}n $\left( \ref{T1 copy(2)}%
\right) $.
\end{theorem}

\section{Problema mixto}

Analizamos ahora las condiciones a imponer para que $U\left( i,j\right) $ en
$\left( \ref{ss}\right) $ verifique el problema mixto $\left( \ref{9}\right)
-\left( \ref{10}\right) .$ Para esto al sustituir $\left( \ref{ss}\right) $
en dichas expresiones, se deben cumplir las siguientes dos ecuaciones:%
\begin{equation}
F\left( i\right) =\sum_{l=1}^{N-1}\left( \text{\ }\hat{E}\hat{E}^{D}P_{l}+%
\hat{E}\hat{E}^{D}Q_{l}\right) v_{l}\left( i\right)  \label{l5.1}
\end{equation}%
\begin{equation}
G\left( i\right) =\frac{\sum_{l=1}^{N-1}\left( \hat{Z}_{0}^{\text{ }}\hat{E}%
\hat{E}^{D}P_{l}+\hat{Z}_{1}^{\text{ }}\hat{E}\hat{E}^{D}Q_{l}\right)
v_{l}\left( i\right) -F\left( i\right) }{k}  \label{l5.2}
\end{equation}%
con lo anterior, seg\'{u}n $\left( \ref{TTT1}\right) ,$ cada componente de $%
F\left( i\right) $ y $G\left( i\right) $ tienen su representaci\'{o}n
discreta respecto a $\left\{ v_{l}\left( i\right) \right\} _{l=1}^{N-1}$,
luego las ecuaciones $\left( \ref{l5.1}\right) $ y $\left( \ref{l5.2}\right)
$ son equivalentes a

\begin{equation}
\hat{E}\hat{E}^{D}P_{l}+\hat{E}\hat{E}^{D}Q_{l}=\frac{\sum_{i=1}^{N-1}v_{l}%
\left( i\right) F\left( i\right) }{\sum_{i=1}^{N-1}v_{l}\left( i\right) ^{2}}%
,  \label{40}
\end{equation}%
\begin{equation}
\hat{Z}_{0}\hat{E}\hat{E}^{D}P_{l}+\hat{Z}_{1}\hat{E}\hat{E}^{D}Q_{l}=\frac{%
\sum_{i=1}^{N-1}v_{l}\left( i\right) \left( kG\left( i\right) +F\left(
i\right) \right) }{\sum_{i=1}^{N-1}v_{l}\left( i\right) ^{2}},  \label{41}
\end{equation}%
que es un sistema matricial cuyas inc\'{o}gnitas son los vectores $P_{l}$ y $%
Q_{l}$ los cuales son soluciones de las ecuaciones%
\begin{equation}
\left( \hat{Z}_{1}-\hat{Z}_{0}^{\text{ }}\right) \hat{E}\hat{E}^{D}P_{l}=%
\frac{\sum_{i=1}^{N-1}v_{l}\left( i\right) \left[ kG\left( i\right) +\left(
\hat{Z}_{1}-I\right) F\left( i\right) \right] }{\sum_{i=1}^{N-1}v_{l}\left(
i\right) ^{2}},  \label{z3}
\end{equation}%
\begin{equation}
\left( \hat{Z}_{1}-\hat{Z}_{0}^{\text{ }}\right) \hat{E}\hat{E}^{D}Q_{l}=%
\frac{\sum_{i=1}^{N-1}v_{l}\left( i\right) \left[ \left( \hat{Z}%
_{0}-I\right) F\left( i\right) -kG\left( i\right) \right] }{%
\sum_{i=1}^{N-1}v_{l}\left( i\right) ^{2}}.  \label{z4}
\end{equation}%
Sin p\'{e}rdida de generalidad, suponiendo que $\sigma \left( D\right) $ no
contiene a $0,$ por el teorema $\left( \ref{AAA}\right) ,$ se tiene que $%
\sigma \left( D\right) =\sigma \left( \hat{E}^{D}\hat{A}\right) -\left\{
0\right\} ,$ adem\'{a}s%
\begin{equation}
\sigma \left( P_{+}\left( D\right) \right) =\sigma \left( \hat{Z}_{0}\right)
-\left\{ 0\right\} ,\text{ }\sigma \left( P_{+}\left( D\right) \right)
=\sigma \left( \hat{Z}_{1}\right) -\left\{ 0\right\}  \label{ez1}
\end{equation}%
o sea
\begin{equation*}
\sigma \left( P_{\pm }\left( D\right) \right) =\left\{ 1+\frac{\rho }{2}d\pm
\sqrt{\left( I+\frac{\rho }{2}d\right) ^{2}-1}:d\in \sigma \left( \hat{E}%
\hat{A}^{D}\right) -\left\{ 0\right\} \right\}
\end{equation*}%
luego%
\begin{equation*}
\sigma \left( P_{-}\left( D\right) -P_{+}\left( D\right) \right) =\left\{ 2%
\sqrt{\left( I+\frac{\rho }{2}d\right) ^{2}-1}:d\in \sigma \left( \hat{E}%
\hat{A}^{D}\right) \right\}
\end{equation*}%
y as\'{\i} $P_{-}\left( D\right) -P_{+}\left( D\right) $ es invertible si se
verifica
\begin{equation}
\rho d\left( 1+\frac{\rho }{4}d\right) \neq 0  \label{ez}
\end{equation}%
Luego tomando en $\left( \ref{ez}\right) $ $\rho $ suficientemente peque\~{n}%
o, utilizando la inversa de $P_{-}\left( D\right) -P_{+}\left( D\right) $ se
pueden emplear los numerales 1 y 3 del teorema $\left( \pageref{INV}\right)
, $ para hallar $\left( \left( \hat{Z}_{1}-\hat{Z}_{0}^{\text{ }}\right)
\hat{E}\hat{E}^{D}\right) ^{G}$, la cual verifica:%
\begin{equation*}
\left( \left( \hat{Z}_{1}-\hat{Z}_{0}^{\text{ }}\right) \hat{E}\hat{E}%
^{D}\right) \left( \left( \hat{Z}_{1}-\hat{Z}_{0}^{\text{ }}\right) \hat{E}%
\hat{E}^{D}\right) ^{G}=\hat{E}\hat{E}^{D}
\end{equation*}%
As\'{\i}, usando el numeral 2 del teorema $\left( \pageref{INV}\right) ,$ se
sigue que $\hat{E}\hat{E}^{D}$ conmuta con $\hat{Z}_{0}$ y $\hat{Z}_{1}^{%
\text{ }}$ y por lo tanto los sistemas $\left( \ref{z3}\right) $ y $\left( %
\ref{z4}\right) $ son consistentes si cumplen la condici\'{o}n%
\begin{equation}
\hat{E}\hat{E}^{D}F\left( i\right) =F\left( i\right) ,\text{ }\hat{E}\hat{E}%
^{D}G\left( i\right) =G\left( i\right)  \label{zz1}
\end{equation}%
teniendo presente lo anterior, por el teorema $\left( \ref{R1.2}\right) $,
las ecuaciones $\left( \ref{z3}\right) $ y $\left( \ref{z4}\right) $ admiten
las siguientes soluciones%
\begin{equation}
P_{l}=\left[ \left( \hat{Z}_{1}-\hat{Z}_{0}^{\text{ }}\right) \hat{E}\hat{E}%
^{D}\right] ^{G}\frac{\sum_{i=1}^{N-1}v_{l}\left( i\right) \left[ kG\left(
i\right) +\left( \hat{Z}_{1}-I\right) F\left( i\right) \right] }{%
\sum_{i=1}^{N-1}v_{l}\left( i\right) ^{2}},  \label{z5}
\end{equation}%
\begin{equation}
Q_{l}=\left[ \left( \hat{Z}_{1}-\hat{Z}_{0}^{\text{ }}\right) \hat{E}\hat{E}%
^{D}\right] ^{G}\frac{\sum_{i=1}^{N-1}v_{l}\left( i\right) \left[ \left(
\hat{Z}_{0}-I\right) F\left( i\right) -kG\left( i\right) \right] }{%
\sum_{i=1}^{N-1}v_{l}\left( i\right) ^{2}}.  \label{z6}
\end{equation}%
Con las condici\'{o}nes
\begin{eqnarray}
&&%
\begin{array}{c}
\left\{ G\left( i\right) ,F\left( i\right) ,\text{ }0<i<N\right\} \subset
KerG\left( \mu ,\beta \right) ,%
\end{array}
\label{z7} \\
&&\text{ \ }%
\begin{array}{c}
G\left( \mu ,\beta \right) \hat{E}^{D}\hat{A}\left( I-G\left( \mu ,\beta
\right) ^{+}G\left( \mu ,\beta \right) \right) =0,%
\end{array}
\notag
\end{eqnarray}%
por el teorema $\left( \ref{T4}\right) ,$ $KerG\left( \mu ,\beta \right) $
es un subespacio invariante de $\left( \hat{E}^{D}\hat{A}\right) .$ Por lo
tanto con $\left( \ref{zz1}\right) $ $P_{l}$ y $Q_{l}$ cumplen%
\begin{equation*}
\hat{E}\hat{E}^{D}P_{l}=P_{l},\text{ }\hat{E}\hat{E}^{D}Q_{l}=Q_{l},\text{ }%
\left\{ P_{l},\text{ }Q_{l},\text{ }0<i<N\right\} \subset KerG\left( \mu
,\beta \right) .
\end{equation*}%
lo que garantiza la consistencia de la condici\'{o}n de contorno
discretizada $\left( \ref{7}\right) -\left( \ref{8}\right) $ junto con el
problema mixto $\left( \ref{9}\right) -\left( \ref{10}\right) .$ Luego, una
soluci\'{o}n de $\left( \ref{6}\right) -\left( \ref{10}\right) $ viene dada
por%
\begin{equation}
U\left( i,j\right) =\sum_{l=1}^{N-1}\left( \hat{Z}_{0}^{\text{\ }j}\hat{E}%
\hat{E}^{D}P_{l}+\hat{Z}_{1}^{\text{\ }j}\hat{E}\hat{E}^{D}Q_{l}\right)
v_{l}\left( i\right) ,  \label{rr1}
\end{equation}%
con $P_{l}$ y $Q_{l}$ determinados por $\left( \ref{z5}\right) $ y $\left( %
\ref{z6}\right) .$

\subsection{Estabilidad}

Para el an\'{a}lisis de la estabilidad de la soluci\'{o}n $\left( \ref{rr1}%
\right) $, suponemos que $F\left( j\right) $ y $G\left( j\right) $ son
acotadas. Por los teoremas $\left( \ref{AE1}\right) $ y $\left( \ref{t0}%
\right) $, de $\left( \ref{ez1}\right) $, el espectro de las matrices $\hat{Z%
}_{0}^{\text{\ }}$ y $\hat{Z}_{1}^{\text{\ }}$ en $\left( \ref{rr1}\right) $
comprende los valores
\begin{equation*}
1+\frac{\rho }{2}d\pm \sqrt{\left( I+\frac{\rho }{2}d\right) ^{2}-1},\text{ }%
\forall d\in \sigma \left( \hat{E}\hat{A}^{D}\right) ,
\end{equation*}%
dado que $\left\vert I+\frac{\rho }{2}d+\sqrt{\left( I+\frac{\rho }{2}%
d\right) ^{2}-I}\right\vert ^{2}=1,$ $\forall d\in \sigma \left( \hat{E}\hat{%
A}^{D}\right) $, considerando que $\rho $ puede verificar%
\begin{equation}
\rho \leq max\left\{ \left\vert \lambda _{l}\right\vert ;\text{ }\lambda
_{l},\text{ }1\leq l\leq N-1\right\} \times r^{2},  \label{eza}
\end{equation}%
por teorema $\left( \ref{sig}\right) $ se obtiene%
\begin{eqnarray}
\left\Vert \hat{Z}_{0}^{\text{\ }}\right\Vert &\leq &\left\Vert \hat{Z}_{1}^{%
\text{\ }}\right\Vert \leq 1+O\left( k\right) ,\left\Vert \hat{Z}_{0}-I^{%
\text{\ }}\right\Vert \leq \left\Vert \hat{Z}_{1}-I^{\text{\ }}\right\Vert
\leq O\left( k\right) \text{ }  \label{r24} \\
&&\text{ \ \ \ \ \ }%
\begin{array}{c}
\left\Vert \left( \hat{Z}_{0}-\hat{Z}_{1}\right) \hat{E}\hat{E}%
^{D}\right\Vert \leq O\left( k^{-1}\right) ,\text{ }k\rightarrow 0.%
\end{array}
\notag
\end{eqnarray}%
para alguna norma matricial $\left\Vert {}\right\Vert $ en $%
\mathbb{C}
^{m\times m}.$ Usando $\left( \ref{z3}\right) ,$ $\left( \ref{z4}\right) $, $%
\left( \ref{r24}\right) $ y el acotamiento de $F(i)$ y $G(i)$, concluimos
que
\begin{equation}
\left\Vert P_{l}\right\Vert =O\left( 1\right) ,\text{ }\left\Vert
Q_{l}\right\Vert =O\left( 1\right) .  \label{r25}
\end{equation}%
De lo anterior$,$ se sigue que $\left\{ U\left( i,j\right) \right\} $ en $%
\left( \ref{rr1}\right) $ permanece acotada, si los n\'{u}meros $\left\Vert
\hat{Z}_{0}^{\text{\ }j}\right\Vert ,$ $\left\Vert \hat{Z}_{1}^{\text{\ }%
j}\right\Vert $ permanecen acotados cuando $j\rightarrow \infty ,$ $%
k\rightarrow 0,$ $1\leq j\leq M,$ $Mk=T.$ N\'{o}tese que de $\left\Vert \hat{%
Z}_{0}^{\text{\ }}\right\Vert \leq 1+O\left( k\right) $ se obtiene que para
alguna constante positiva $S,$ $\left\Vert \hat{Z}_{0}^{\text{\ }%
}\right\Vert \leq 1+kS,$ entonces, para $1\leq j\leq M$ se calcula%
\begin{equation}
\left\Vert \hat{Z}_{0}^{\text{\ }j}\right\Vert \leq \left\Vert \hat{Z}_{0}^{%
\text{\ }}\right\Vert ^{j}\leq \left( 1+O\left( k\right) \right) ^{j}\leq
\left( 1+O\left( k\right) \right) ^{M}\leq e^{MO\left( k\right) }\leq
e^{MkS}=e^{TS}  \label{r26}
\end{equation}%
lo mismo ocurre para $\left\Vert \hat{Z}_{1}^{\text{\ }j}\right\Vert .$ As%
\'{\i}, considerando $\left( \ref{r25}\right) $, $\left( \ref{r26}\right) $
y $L=\max \left\{ \left\Vert v_{l}\right\Vert :1\leq l\leq N-1\right\} ,$ se
sigue que la soluci\'{o}n definida por $\left( \ref{rr1}\right) ,$ $\left( %
\ref{z5}\right) $ y $\left( \ref{z6}\right) $ es estable, es decir:%
\begin{eqnarray*}
\left\Vert U\left( i,j\right) \right\Vert &=&O\left( 1\right) ,\text{ }%
k\rightarrow 0,\text{ }h=\frac{1}{N}\text{ paso temporal fijo} \\
1 &\leq &i\leq N-1\text{ \ }j\rightarrow \infty ,\text{ }con\text{ }Mk=T
\end{eqnarray*}%
En resumen, el siguiente resultado ha sido establecido.

\begin{theorem}
\label{Main}Con las mismas consideraciones del teorema anterior, sean $%
G\left( i\right) $ y $F\left( i\right) $ acotadas que satisfacen $\left(
\text{\ref{zz1}}\right) $ y $\left( \text{\ref{z7}}\right) $ con $\rho $
verificando $\left( \ref{ez}\right) $ y $\left( \ref{eza}\right) $;
entonces\ una soluci\'{o}n estable del problema discreto $\left( \ref{6}%
\right) -\left( \ref{10}\right) $ esta dada, por%
\begin{equation}
U\left( i,j\right) =\sum_{l=1}^{N-1}\left( \hat{Z}_{0}^{\text{\ }j}\hat{E}%
\hat{E}^{D}P_{l}+\hat{Z}_{1}^{\text{\ }j}\hat{E}\hat{E}^{D}Q_{l}\right)
v_{l}\left( i\right)  \label{zzz1}
\end{equation}%
donde $P_{l},$ $Q_{l}$ son como lo indican las ecuaciones $\left( \ref{z5}%
\right) $ y $\left( \ref{z6}\right) .$
\end{theorem}

\subsection{Ejemplo}

Consideremos las matrices en $%
\mathbb{C}
^{3\times 3}$%
\begin{eqnarray*}
&&%
\begin{array}{c}
\text{ \ \ }E=\left(
\begin{array}{ccc}
\epsilon & 0 & \delta \\
0 & \gamma & 0 \\
0 & 0 & 0%
\end{array}%
\right) ,\text{ }A=\left(
\begin{array}{ccc}
0 & 0 & 0 \\
0 & \delta & 0 \\
0 & 0 & \sigma%
\end{array}%
\right) ,\text{ }\left\{ \epsilon ,\text{ }\gamma ,\text{ }\sigma \right\}
\text{ }\in
\mathbb{C}
-\left\{ 0\right\} \text{.}%
\end{array}%
\text{ } \\
&&\text{ }%
\begin{array}{c}
A_{1}=\left(
\begin{array}{ccc}
a_{11} & a_{12} & b_{1} \\
a_{21} & a_{22} & b_{2} \\
a_{31} & a_{32} & b_{3}%
\end{array}%
\right) ,\text{ }A_{2}=\left(
\begin{array}{ccc}
\mu a_{11} & \mu a_{12} & c_{1} \\
\mu a_{21} & \mu a_{22} & c_{2} \\
\mu a_{31} & \mu a_{32} & c_{3}%
\end{array}%
\right) ,\text{ }B_{2}=\eta B_{1}%
\end{array}%
\end{eqnarray*}%
y las funciones vectoriales discretas%
\begin{equation*}
F(i)=\left(
\begin{array}{c}
f_{1}\left( i\right) \\
f_{2}\left( i\right) \\
0%
\end{array}%
\right) ,\text{ }G(i)=\left(
\begin{array}{c}
g_{1}\left( i\right) \\
g_{2}\left( i\right) \\
0%
\end{array}%
\right) ,
\end{equation*}%
escogiendo $\lambda \neq 0$ tal que $\lambda \gamma +\delta \neq 0$, se
calcula%
\begin{eqnarray*}
&&%
\begin{array}{c}
\text{ \ \ \ \ \ \ \ \ \ \ \ \ \ \ \ \ \ \ \ \ \ \ \ \ \ \ }\left( \lambda
E+A\right) ^{-1}=\left(
\begin{array}{ccc}
1/\lambda \epsilon & 0 & -\delta /\epsilon \sigma \\
0 & 1/\left( \lambda \gamma +\delta \right) & 0 \\
0 & 0 & 1/\sigma%
\end{array}%
\right) ,\text{ }%
\end{array}
\\
&&%
\begin{array}{c}
\text{ \ \ \ \ \ \ \ \ \ \ \ }\hat{E}=\left(
\begin{array}{ccc}
1/\lambda & 0 & \delta /\lambda \epsilon \\
0 & \gamma /\left( \lambda \gamma +\delta \right) & 0 \\
0 & 0 & 0%
\end{array}%
\right) ,\text{ }\hat{A}=\left(
\begin{array}{ccc}
0 & 0 & -\delta /\epsilon \\
0 & \delta /\left( \lambda \gamma +\delta \right) & 0 \\
0 & 0 & 1%
\end{array}%
\right) ,%
\end{array}
\\
&&%
\begin{array}{c}
\hat{E}^{D}=\left(
\begin{array}{ccc}
\lambda & 0 & \delta \lambda /\epsilon \\
0 & \left( \lambda \gamma +\delta \right) /\gamma & 0 \\
0 & 0 & 0%
\end{array}%
\right) ,\text{ }\hat{E}\hat{E}^{D}=\left(
\begin{array}{ccc}
1 & 0 & \delta /\epsilon \\
0 & 1 & 0 \\
0 & 0 & 0%
\end{array}%
\right) ,\text{ }\hat{E}^{D}\hat{A}=\left(
\begin{array}{ccc}
0 & 0 & 0 \\
0 & \delta /\gamma & 0 \\
0 & 0 & 0%
\end{array}%
\right) ,%
\end{array}%
\end{eqnarray*}%
si en $\left( \ref{32a}\right) -\left( \ref{34a}\right) $ escogemos$,$ $%
\alpha =\mu $ y $\beta =\eta $, entonces%
\begin{equation*}
\text{ }G\left( \alpha ,\beta \right) =\left(
\begin{array}{cccccc}
0 & 0 & 0 & 0 & 0 & 0 \\
0 & 0 & 0 & 0 & 0 & 0 \\
d_{1} & d_{2} & d_{3} & 0 & 0 & 0%
\end{array}%
\right) ^{T},\text{ }G\left( \alpha ,\beta \right) ^{+}=1/\left(
d_{1}^{2}+d_{1}^{2}+d_{1}^{2}\right) \left(
\begin{array}{cccccc}
0 & 0 & 0 & 0 & 0 & 0 \\
0 & 0 & 0 & 0 & 0 & 0 \\
d_{1} & d_{2} & d_{3} & 0 & 0 & 0%
\end{array}%
\right) ,
\end{equation*}%
donde $d_{i}=\mu b_{i}-c_{i},$ $i=1,$ $2,3$ y se comprueba adem\'{a}s que:%
\begin{eqnarray*}
&&%
\begin{array}{c}
\sigma \left( \hat{E}^{D}\text{\ }\hat{A}\right) =\left\{ 0,\text{ }1,\text{
}3\right\} ,ran\left( G\left( 2,1\right) \right) <3,\text{ }G\left( \alpha
,\beta \right) F(i)=0,\text{ }G\left( \alpha ,\beta \right) G(i)=0,%
\end{array}
\\
&&\text{ }%
\begin{array}{c}
\hat{E}\hat{E}^{D}F(i)=F(i),\text{ \ }\hat{E}\hat{E}^{D}G(i)=G(i),\text{ }%
G\left( \mu ,\beta \right) \hat{E}^{D}\hat{A}\left( I-G\left( \mu ,\beta
\right) ^{+}G\left( \mu ,\beta \right) \right) =0.%
\end{array}%
\end{eqnarray*}%
Luego, si $\rho $ es suficientemente peque\~{n}o, las hip\'{o}tesis del
teorema $\left( \ref{Main}\right) $ pueden ser verificadas y por lo tanto
una soluci\'{o}n de $\left( \ref{6}\right) -\left( \ref{10}\right) $ es de
la forma%
\begin{equation*}
U\left( i,j\right) =\sum_{l=1}^{N-1}\left( \hat{Z}_{0}^{\text{\ }j}\hat{E}%
\hat{E}^{D}P_{l}+\hat{Z}_{1}^{\text{\ }j}\hat{E}\hat{E}^{D}Q_{l}\right)
v_{l}\left( i\right)
\end{equation*}%
deonde $P_{l}$ y $Q_{l}$ son como lo indican las ecuaciones $\left( \ref{z5}%
\right) $ y $\left( \ref{z6}\right) $ y $v_{l}\left( i\right) $ son las
eigenfunciones descritas en la proposici\'{o}n $\left( \ref{T1 copy(2)}%
\right) .$

\subsection{Observaci\'{o}n}

\begin{enumerate}
\item El procedimiento para el hallazgo de la soluci\'{o}n del problema $%
\left( \ref{1}\right) -\left( \ref{4}\right) $ aqu\'{\i} estudiado, incluye
el caso de su versi\'{o}n no singular. En efecto, sin perder generalidad,
supongamos que en la ecuaci\'{o}n $\left( \ref{1}\right) ,$ $E=I.$ Al
considerar la descomposici\'{o}n can\'{o}nica de $A=T\left(
\begin{array}{cc}
K & 0 \\
0 & N%
\end{array}%
\right) T^{-1}$, la matriz $\lambda I+A$ obtiene la forma $\lambda
I+A=T\left(
\begin{array}{cc}
\lambda I+K & 0 \\
0 & \lambda I+N%
\end{array}%
\right) T^{-1}.$ Escogiendo $\lambda $ suficientemente grande, esto har\'{a}
que la diagonal de los bloques $\lambda I+K$ y $\lambda I+N$ no se anule y as%
\'{\i}, la condici\'{o}n: que existe $\gamma \in
\mathbb{C}
$ tal que $(\gamma E+A)$ es no singular, es verificada, y nuestro
procedimiento puede ser aplicado. N\'{o}tese adem\'{a}s que $\hat{E}^{D}\hat{%
E}=I$ y $\hat{E}^{D}\hat{A}=A,$ por lo tanto la condici\'{o}n $\left( \ref%
{zz1}\right) $ puede omitirse$.$

\item Si la condici\'{o}n $\left( \ref{5}\right) $ del teorema $\ref{Main}$
no se cumple, entonces de acuerdo a la proposici\'{o}n $\ref{AAA},$ la ecuaci%
\'{o}n$\ \left( \ref{16}\right) $ puede ser reducida a $\left( \ref{25}%
\right) ;$ como la soluci\'{o}n de esta es la soluci\'{o}n trivial, se
concluye que as\'{\i} tambi\'{e}n lo ser\'{a} la soluci\'{o}n encontrada
para $\left( \ref{6}\right) -\left( \ref{10}\right) .$
\end{enumerate}

\subsection{Conclusion}

En este trabajo se ha demostrado la existencia de soluciones num\'{e}ricas
estables para una gama de sistemas de ecuaciones en derivadas parciales
singulares de la ecuaci\'{o}n de onda homog\'{e}nea como la descrita en $%
\left( \ref{1}\right) -\left( \ref{4}\right) $. En nuestra metodolog\'{\i}a,
con respecto a los eigenvalores del problema de Sturm Liouville discreto
asociados a $\left( \ref{32a}\right) -\left( \ref{34a}\right) $, hemos
asumido el teorema espectral para matrices sim\'{e}tricas en vez de aplicar m%
\'{e}todos num\'{e}ricos para establecer su existencia. Dicho enfoque nos
permiti\'{o} enfrentar con otra estrategia el problema de encontrar tal
soluci\'{o}n a nuestro sistema de inter\'{e}s para el cual la arbitrariedad
de los par\'{a}metros $\alpha $ y $\beta ,$ relacionados al PSD $\left( \ref%
{32a}\right) -\left( \ref{34a}\right) $ puede ser usada con el fin de que
las hip\'{o}tesis del teorema $\left( \ref{Main}\right) $ se verifiquen.


\begin{thebibliography}{99}
\bibitem{Argawal} R. P. Agarwual. Difference Equations and Inequalities.
Theory, Methods and Applications. Marcel Dekker, 1992.%
\addcontentsline{toc}{chapter}{Bibliografía}

\bibitem{P1} S.L . Campbell and C.D. Meyer, Jr., Generalized inverses of
Linear Transformations. Pitman. London, 1979.

\bibitem{camacho} J. Camacho, E. Defez, L. J\'{o}dar and J. V. Romero.
Discrete numerical solution of coupled mixed hyperbolic problems, Computers
\& Mathematics with Applications 46, 8-9 (2003), 1183-1193.

\bibitem{S10} M. C. Casab\'{a}na, L. Jodar, G.A. Ossand\'{o}n. Conditional
uniform time stable numerical solutions of coupled hyperbolic systems.
Applied Mathematics Letters 20 (2007) 13--16.

\bibitem{S3} P. K. Das, Optical Signal Processing, Springer, New York,
(1991).

\bibitem{s2} A. Das, Sisir K. Das, Microwave Engineering, Mc Graw-Hill, West
Patel Nagar (2000).

\bibitem{Dunford} N. Dunford and J. Schwartz. Linear Operators. Part I.
Interscience, New York, 1988.

\bibitem{Analitic Camacho} L. J\'{o}dar, E. Navarro, J. Camacho.
Analytic-numerical solutions with a priori error bounds for a class of
strongly coupled mixed partial differential systems. Journal of
Computational and Applied Mathematics 104, 2 (1999), 123-143.

\bibitem{25a} L. J\'{o}dar, E. Navarro, A. E. Posso, M. C. Casab\'{a}n.
Constructive solution of strongly coupled continuous hyperbolic mixed
problems. Applied Numerical Mathematics 47, 3-4 (2003), 477 - 492.

\bibitem{PON} http://paginasweb.univalle.edu.co/\symbol{126}%
ccm2009/publico/ponencia.php?idioma=ES\&id=838.

\bibitem{P34} G. Ossandon. Construcci\'{o}n de soluciones num\'{e}ricas
estables de sistemas en derivadas parciales fuertemente acoplados mediante m%
\'{e}todos semi-impl\'{\i}citos. Tesis doctoral, Universidad Polit\'{e}cnica
de Valencia, 2004.

\bibitem{P340a} Robert Piziak, Patrick L. Odell. Matrix theory: from
generalized inverses to Jordan form. Chapman \& Hall/CRC, Virginia, 2007.

\bibitem{S1} A. C Metaxas and R. J. Meredith, Indutrial Microwave Heating,
Peter Peregrinus, London, (1983).

\bibitem{P10a} G. D. Smith. Numerical Solution of Partial Differential
Equations: finite difference methods, Oxford University Press, 1985.

\bibitem{S11} M. A Slawinsky, Seismic Waves and Rays in Elastic Media
(Handbook of Geophysical Exploration: Seismic Exploration), Elsevier
Science, Oxford, 2003.

\bibitem{P37a} E. E. Villa Chica. Soluciones discretas estables de problemas
mixtos singulares para sistemas hiperb\'{o}licos fuertemente acoplados\emph{%
. }Tesis de maestr\'{\i}a, Universidad de Antioquia, Medell\'{\i}n, 2006.

\bibitem{S7a} A.T. Winfree, When Times Breaks Dowmn, Princeton Univ. Press,
Princeton, 1987.
\end{thebibliography}
\end{document}